# Optimization of Fuzzy Random Portfolio selection by Implementation of Harmony Search Algorithm

Mir Ehsan Hesam Sadati[#1], Ali Doniavi[*2]

[#] *Department of Industrial Engineering*
*Urmia University, Urmia, Iran*

*Abstract*— **This study first reviews fuzzy random Portfolio selection theory and describes the concept of portfolio optimization model as a useful instrument for helping finance practitioners and researchers. Second, this paper specifically aims at applying possibility-based models for transforming the fuzzy random variables to the linear programming. The harmony search algorithm approaches to resolve the portfolio selection problem with the objective of return maximization is applied. We provide a numerical example to illustrate the proposed model. The results show that the evolutionary method of this paper with harmony search algorithm, can consistently handle the practical portfolio selection problem.**

## I. INTRODUCTION

The mean-variance model originally introduced by Markowitz [1] plays an important role in the development of modern portfolio selection theory. Portfolio optimization (PO) consists of the portfolio selection problem in which we want to find the optimum way of investing a particular amount of money in a given set of securities or assets [2]. In many industries, there are many decision problems; i.e., scheduling problem, logistics. In these problems, it is important to predict future total returns and to decide an optimal asset allocation maximizing total profits under some constraints. We call such industrial assets allocation problems portfolio selection problems.

Markowitz formulated mean-variance models mathematically in two ways: minimizing variance for a given expected value, or maximizing expected value for a given variance. When selecting portfolio, an investor may encounter with both fuzziness and randomness. In fact, for an investor, the fuzziness and randomness of security returns are often mixed up with each other. In such situations, we may employ fuzzy random theory [3] to deal with this uncertainty of fuzziness and randomness. A Fuzzy random variable was first introduced by Kwakernaak [4], and its mathematical basis was constructed by Puri and Ralescu [5]. In this paper, the asset return in portfolio selection problem are fuzzy random variables and we use the concept of possibility-based model to develop a solution method for the fuzzy random portfolio optimization problem.

Geem et al. [6] proposed a new meta-heuristic algorithm, so-called the Harmony Search (HS) algorithm. Although it is comparatively a new meta-heuristic algorithm, in various applications, it has been proven to be a robust and efficient tool.

The rest of the paper is organized as follows: Section 2 includes basic concepts on fuzzy and fuzzy random theory. In Section 3, the problem formulation is presented. Section 4, explains the Harmony Search meta-heuristic algorithm in details. In section 5, numerical example is solved to illustrate the proposed model. Finally conclusion and future work will be presented in section 6.

## II. BASIC CONCEPTS

The term fuzzy random variable was coined by Kwakernaak [4], who introduced FRVs as "random variables whose values are not real, but fuzzy numbers," and conceptualized a FRV as a vague perception of a crisp but unobservable RV, and its mathematical basis was constructed by Puri and Ralescu [5].

The concept of fuzzy random variable was introduced as an analogous notion to random variable in order to extend statistical analysis to situations when the outcomes of some random experiment are fuzzy sets. In general, fuzzy random variables can be defined in an n dimensional Euclidian space $R^n$. We present the definition of a fuzzy random variable in a single dimensional Euclidian space $R$.

*Definition 1*

LR fuzzy number $\tilde{A}$ is defined by following membership function:

$$\tilde{A}(X) = \begin{cases} L\left(\dfrac{A^0 - x}{\beta}\right) & if \quad A^0 - \beta \le x \le A^0 \\ 1 & if \quad A^0 \le x \le A^1 \\ R\left(\dfrac{x - A^1}{\gamma}\right) & if \quad A^1 \le x \le A^1 + \gamma \end{cases} \quad (1)$$

where $[A^0, A^1]$ show the peak of fuzzy number $\tilde{A}$ and $\beta, \gamma$ represent the left and right spread respectively; $L, R = [0,1] \to [0,1]$ with $L(0) = R(0) = 1$ and $L(1) = R(1) = 0$ are strictly decreasing, continuous functions. A possible representation of a LR fuzzy number is $\tilde{A} = (A^0, A^1, \beta, \gamma)_{LR}$.

*Definition 2*

Let $(\Omega, A, P)$ be a probability space, where $\Omega$ is a sample space, A is a σ-field and P is a probability measure. Let $F_N$ be the set of all fuzzy numbers and B a Borel σ-field of R. Then a map $\tilde{\tilde{Z}} : \Omega \to F$ is called a fuzzy random variable if it holds that

$$\{(\omega, \tau) \in \Omega \times R \mid \tau \in \tilde{\tilde{Z}}_\alpha(\omega)\} \in A \times B, \forall \alpha \in [0,1] \quad (2)$$

where

$$\tilde{\tilde{Z}}_\alpha(\omega) = \left[\tilde{\tilde{Z}}_\alpha^-(\omega), \tilde{\tilde{Z}}_\alpha^+(\omega)\right] = \left\{\tau \in R \mid \mu_{\tilde{Z}(\omega)}(\tau) \geq \alpha\right\} \quad (3)$$

is an α-level set of the fuzzy number $\tilde{Z}(\omega)$ for $\omega \in \Omega$.

Intuitively, fuzzy random variables are considered to be random variables whose realized values are not real values but fuzzy numbers or fuzzy sets.

III. PROBLEM FORMULATION

In the following problem, like the mean-variance model introduced by Markowitz [1] is called Fuzzy Random portfolio selection problem, the return rate of assets are fuzzy random variables:

Problem 1

$$Max \ \tilde{\tilde{Z}} = \sum_{j=1}^n \tilde{\tilde{R}}_j x_j \quad (4)$$

$$s.t. \quad \sum_{j=1}^n x_j = M_0, \quad (5)$$

$$\sum_{j=1}^n \tilde{\tilde{R}}_j x_j \geq \tilde{\tilde{R}}_0 \quad (6)$$

$$0 \leq x_j \leq U_j ; \quad j = 1, 2, ..., n. \quad (7)$$

The parameters and variables are define as follow, for j=1, 2,…, n:

$\tilde{\tilde{R}}_j = (R_j^0, R_j^1, \beta_j, \gamma_j)_{LR}$ represents fuzzy random variables whose observed value for each $\omega \in \Omega$ is fuzzy number $\tilde{R}_j(\omega) = (R_j^0(\omega), R_j^1(\omega), \beta_j, \gamma_j)_{LR}$.

$(\overline{R}_j^0, \overline{R}_j^1) = (R_j^0 + \overline{t}R_j^2, R_j^1 + \overline{t}R_j^2)$ is a random vector in which $\overline{t}$ is a random variable with cumulative distribution function $T$.

$n$: The number of assets for possible investment

$M_0$: Available total fund

$\tilde{\tilde{R}}_j$: The rate of return of asset j (per period)

$\tilde{\tilde{R}}_0$: The return in dollars

$x_j$: Decision variables which represent the dollar amount of fund invested in asset j

$U_j$: The upper bound of investment in asset j.

*Possibility-based Model*

By Zadeh's extension principle for objective function in problem 1, its membership function is given as follows for each $\omega \in \Omega$:

$$\mu_{\tilde{Z}(\omega)}(t) = \begin{cases} L\left(\dfrac{Z^0(\omega) - t}{\beta}\right) & if \quad t \leq Z^0(\omega) \\ 1 & if \quad Z^0(\omega) \leq t \leq Z^1(\omega) \\ R\left(\dfrac{t - Z^1(\omega)}{\gamma}\right) & if \quad otherwise \end{cases} \quad (8)$$

where $\tilde{Z}(\omega) = (Z^0(\omega), Z^1(\omega), \beta, \gamma)_{LR}, Z^0(\omega) = \sum_{j=1}^n R_j^0(\omega) x_j$, and $Z^1(\omega) = \sum_{j=1}^n R_j^1(\omega) x_j$.

The degree of possibility $\pi(\tilde{Z}(\omega) \geq f)$ under the possibility distribution $\mu_{\tilde{Z}(\omega)}(t)$ is given as follows:

$$\pi(\tilde{Z}(\omega) \geq f) = \sup_{y_1, y_2}\left\{\min\left\{\mu_{\tilde{Z}(\omega)}(y_1), \mu_f(y_2)\right\} \mid y_1 \geq y_2\right\} \geq \eta. \quad (9)$$

The possibility degree of fuzzy constraint $\left(\sum_{j=1}^n \tilde{R}_j(\omega) x_j \geq \tilde{R}_0(\omega)\right)$ under the possibility distributions is defined as follows:

$$\pi\left(\sum_{j=1}^n \tilde{R}_j(\omega) x_j \geq \tilde{R}_0(\omega)\right) =$$
$$\sup_{y_1, y_2}\left\{\min\left\{\mu_{\sum_{j=1}^n \tilde{R}_j(\omega) x_j}(y_1), \mu_{\tilde{R}_0(\omega)}(y_2)\right\} \mid y_1 \geq y_2\right\} \quad (10)$$

We maximize the degree of possibility $\pi(\tilde{Z}(\omega) \geq f)$ and the degree of possibility $\pi\left(\sum_{j=1}^{n} \tilde{R}_j(\omega) x_j \geq \tilde{R}_0(\omega)\right)$, so our portfolio selection model in Problem 1 comes by the following model:

Problem 2

$$Max \ f \qquad (11)$$

$$s.t. \quad \Pr\{\omega \mid \pi(\tilde{Z}(\omega) \geq f) \geq \eta\} \geq \lambda, \qquad (12)$$

$$\sum_{j=1}^{n} x_j = M_0, \qquad (13)$$

$$\Pr\left\{\omega \mid \pi\left(\sum_{j=1}^{n} \tilde{R}_j(\omega) x_j \geq \tilde{R}_0(\omega)\right) \geq \eta\right\} \geq \lambda, \qquad (14)$$

$$0 \leq x_j \leq U_j \ ; \quad j = 1, 2, ..., n. \qquad (15)$$

where $\lambda$ is a predetermined probability level and $\eta$ is a predetermined possibility level. A feasible solution of portfolio selection problem is called a possibility solution. In order to transform the above model to a linear programming model, we need to reformulate (12) and (14). Consider the following theorem:

*Theorem 1: [7]*

For any decision variable, it holds that:

$$\begin{cases} 1) \Pr\{\omega \mid \pi(\tilde{Z}(\omega) \geq f) \geq \eta\} \geq \lambda \Leftrightarrow \\ \sum_{j=1}^{n}\left(R_j^1 + T^*(1-\lambda) R_j^2\right) x_j + R^*(\eta) \sum_{j=1}^{n} \gamma_j x_j \geq f \\ 2) \Pr\left\{\omega \mid \pi\left(\sum_{j=1}^{n} \tilde{R}_j(\omega) x_j \geq \tilde{R}_0(\omega)\right) \geq \eta\right\} \geq \lambda \Leftrightarrow \\ \sum_{j=1}^{n}\left(R_j^1 + T^*(1-\lambda) R_j^2\right) x_j + R^*(\eta) \sum_{j=1}^{n} \gamma_j x_j \geq R_0^0 + T^*(1-\lambda) R_0^2 - \beta_0 L^*(\eta) \end{cases}$$

where $T^*$, $L^*$ and $R^*$ are pseudo inverse functions defined as:

$$T^*(\lambda) = \inf\{t \mid T(t) \geq \lambda\},$$

$$L^*(\lambda) = \sup\{t \mid L(t) \geq \lambda\},$$

$$R^*(\lambda) = \sup\{t \mid R(t) \geq \lambda\}.$$

Now the optimal solution of Problem 2 is equal to the following linear fractional programming problem:

Problem 3

$$Max \ (P) = \sum_{j=1}^{n}\left(R_j^1 + T^*(1-\lambda) R_j^2\right) x_j + R^*(\eta) \sum_{j=1}^{n} \gamma_j x_j \qquad (16)$$

$$s.t. \quad \sum_{j=1}^{n} x_j = M_0, \qquad (17)$$

$$\sum_{j=1}^{n}\left(R_j^1 + T^*(1-\lambda) R_j^2\right) x_j + R^*(\eta) \sum_{j=1}^{n} \gamma_j x_j \geq R_0^0 + T^*(1-\lambda) R_0^2 - \beta_0 L^*(\eta), \qquad (18)$$

$$0 \leq x_j \leq U_j \ ; \quad j = 1, 2, ..., n. \qquad (19)$$

IV. HARMONY SEARCH ALGORITHM

Harmony search is a music-based meta-heuristic optimization algorithm [6]. It was inspired by the observation that the aim of music is to search for a perfect state of harmony. This harmony in music is analogous to find the optimality in an optimization process. The search process in optimization can be compared to a jazz musician's improvisation process. On the one hand, the perfectly pleasing harmony is determined by the audio aesthetic standard. A musician always intends to produce a piece of music with perfect harmony. On the other hand, an optimal solution to an optimization problem should be the best solution available to the problem under the given objectives and limited by constraints [8].

Geem et al. [9] formalized these three options into quantitative optimization process in 2009, and the three corresponding components become: usage of harmony memory, pitch adjusting, and randomization. Similarly, when each decision variable picks a value, there are three options: (1) to pick any value from the memory; (2) to pick a value adjacent to any value in the memory; (3) to pick a random value from the domain of all possible values.

Having explained the three main components of the HS algorithm: harmony memory (HM), harmony memory consideration rate (HMCR) and pitch adjustment rate (PAR), the following subsections explain each step that comprises the HS algorithm [10]:

*A. Problem formulation*

The HS algorithm was initially conceived for solving optimization problems where a single objective and several constraints are considered.

*B. Parameter configuration*

Furthermore, besides the two parameters already mentioned, HMCR and PAR, the HS algorithm has other parameters such as: harmony memory size (HMS), maximum amount of improvisations or iterations (Maximum Improvisations, MI) and pitch range variability (Fret Width, FW) that operate altogether with PAR in pitch adjustment.

*C. Memory initialization*

After the problem has been formulated and the parameters properly configured, a random configuration process is performed on the memory. The HS algorithm initially improvises several solutions randomly. The number of solutions must be at least equal to HMS. Then, the best HMS solutions are selected.

*D. Improvisation*

As it was mentioned, there are three options among which the HS algorithm may choose when performing an improvisation:

**(i) Random selection:** When HS determines the $x_i^{new}$ value for a new solution $x^{new} = \left[x_1^{new}, x_2^{new}, ..., x_n^{new}\right]$, it randomly chooses a value from the range of all possible values $\{x_i(1), x_i(2), ..., x_i(K_i)\}$ or $x_i^l \leq x_i \leq x_i^u$ with a probability of (1 − HMCR).

**(ii) Memory consideration:** When HS determines the value of $x_i^{new}$, it randomly chooses the $x_i^j$ value from HM (j = 1, 2,..., HMS) with a probability equal to HMCR

**(iii) Pitch adjustment:** After the value of $x_i^{new}$ has been randomly chosen from HM in the process previously described, it may be adjusted to neighboring values adding or subtracting a given amount, with probability PAR. For discrete variables, if $x_i(k) = x_i^{new}$, the pitch adjustment is $x_i(k+m)$, where $m \in \{-1, 1\}$. For continuous variables, the pitch adjustment is $x_i^{new} + \Delta$, where $\Delta = U(-1,1) * FW(i)$.

The three components of the HS algorithm described in the above section can easily be implemented using any programming language, though it should be straightforward to carry out simulations with visualization using Matlab.

*E. Memory update*

If the new solution $x^{new}$ is better than the worst solution in HM in terms of the objective function value, the new solution is included in HM and the worst is discarded.

*F. Termination*

If the HS algorithm meets the stopping criterion (for instance, has reached the maximum amount of iterations or the maximum execution time), the process is terminated.

V. AN EXAMPLE

In this section, an example is given to illustrate the proposed harmony search algorithm for portfolio optimization selection by possibility-based model. Let us consider 5 securities whose returns are fuzzy random variables and their values are given in Table 1. $\bar{t}$ is a normal random variable whose mean 0 and variance 1.

The upper bound of investment amount in each stock is set to no more than 60 units of the total available fund. Given a total allocation budget of 200 units and annual return which is fuzzy random variable is shown as $\tilde{\bar{R}} = M_0 \tilde{\bar{r}}_0$ where $\tilde{\bar{r}}_0 = (1 + 0.3\bar{t}, 1 + 0.3\bar{t}, 0.3, 0.3)$. Now we want to know what is the optimal solution for our portfolio selection problem for the different levels of probability and possibility {0.1, 0.4, 0.7, 0.9}.

We apply the harmony search algorithm and possibility-based model based on theorems 1 to obtained optimum solution. We have used 6 harmonics, the harmony accepting rate HMCR=0.9, and the pitch adjusting rate PAR=0.5. We can see that the pitch adjustment is more intensive in local regions (two thin strips), this is probably another reason why the harmony search is more efficient than other algorithms. The best estimate solutions are obtained after 10,000 iterations using the Matlab program. All the results are collected in Table 2.

TABLE I
PARAMETERS OF THE EXAMPLE

| $j$ | 0 | 1 | 2 | 3 | 4 | 5 |
|---|---|---|---|---|---|---|
| $R_j^0$ | 250 | 1.2 | 1.25 | 1.35 | 1.25 | 1.4 |
| $R_j^1$ | 250 | 1.35 | 1.3 | 1.45 | 1.35 | 1.5 |
| $R_j^2$ | 50 | 0.5 | 0.6 | 0.55 | 0.4 | 0.5 |
| $\beta_j$ | 40 | 0.15 | 0.1 | 0.2 | 0.15 | 0.2 |
| $\gamma_j$ | 40 | 0.15 | 0.1 | 0.2 | 0.15 | 0.2 |

TABLE III
NUMERICAL RESULTS (POSSIBILITY-BASED MODEL)

| $\lambda, \eta$ | 0.1 | 0.4 | 0.7 | 0.9 |
|---|---|---|---|---|
| $x_1^*$ | 20 | 0 | 0 | 20 |
| $x_2^*$ | 60 | 60 | 20 | 0 |
| $x_3^*$ | 60 | 60 | 60 | 60 |
| $x_4^*$ | 0 | 20 | 60 | 60 |
| $x_5^*$ | 60 | 60 | 60 | 60 |
| OFV[ap] | 451.22 | 331.85 | 244.39 | 164.44 |

Clearly, the greater the $\lambda, \eta$ value, the greater the level of possibility and the lower the objective function value is.

## VI. CONCLUSION

Portfolio optimization has been one of the important fields of research in economics and finance. Since the prospective returns of assets used for portfolio optimization problem are forecasted values, considerable uncertainty is involved. In this paper, Markowitz's mean-variance idea was extended to portfolio selection by possibility-based model. This paper proposed a solution method for portfolio selection model whose parameters were fuzzy random variables. The idea was based on possibility-based model and the way of finding the optimum solution was harmony search algorithm. The pre-defined parameters for a harmony search were altered, and a variable HMCR, variable PAR and variable bandwidth were utilized. The harmony search was found to be an efficient and robust algorithm. For future research, we will apply the other methods for fuzzy random portfolio selection model and improve our portfolio selection problem.